\documentclass[10pt]{article}
\usepackage{amssymb}
\usepackage{enumerate}
\usepackage{graphicx}
\usepackage{amsthm}
\usepackage{color}
\usepackage{epstopdf}

\def\R{{\rm I\! R}}

\def\dr{{\partial \over \partial r}}

\newtheorem{theorem}{Theorem}
\newtheorem{corollary}{Corollary}

\title{On the uniqueness of limit cycles for Li\'enard equation: the legacy of G. Sansone}

\author{M. Sabatini, G. Villari  }
\begin{document}
\maketitle

\section{Introduction}The aim of this paper is to present some classical and more recent results concerning
the uniqueness of limit cycles for the Lienard equation
\begin{equation} 
\ddot x  + f(x) \dot x + g(x)  = 0 .
\end{equation} 
Such an equation is without any doubt among the most investigated ones in the qualitative theory
of ordinary differential equations, for its applications in mechanics and electric circuits theory.
It may be considered the starting point in the theory of limit cycles, becoming sort of a benchmark for new methods.
Its relationships with different classes of systems, as quadratic ones (which can be reduced to suitable Li\'enard systems) made it even more relevant in the study of planar systems.

We do not intend to make a survey on this widely investigated area.
Our attempt is rather directed to focus on the results of Giovanni Sansone and to
show how in some sense those results were giving a crucial contribute for the development
of this field.  After such equation was proposed by Lienard in 1928 \cite{L}, who proved the first uniqueness theorem for limit cycles, it was necessary to wait until 1942 to find another significant result  in
 this direction, due to Levinson and Smith  \cite{LS}. 
 But it was only with the contribution of Sansone, starting from his celebrated paper of 1949 \cite{San1},
that the relevance of this problem was recognized, not only from the Italian school. We just recall
 the famous paper of Filippov \cite{F}, which was devoted to the problem of existence of limit cycles
 and therefore will be not examined in this context. For the enormous contribution of the Chinese School, which for long time was not known out of China for language problems, we refer to the monograph of Ye Yan Qian \lq\lq Theory of limit cycles" \cite{YYQ}.
 Finally, a survey updated until the 70's is due to Staude \cite{Sta}. Coming back to the work of Sansone, we believe that our, in some sense historical, approach will be of some interest, because it shows how deep results may come out from simple geometrical ideas in a far to be easy field of research as is the problem of uniqueness for periodic solutions. In particular, we would like to describe the three approaches applied by Sansone in such a field:
the energy integral, the divergence integral and, implicitly, the rotation of the vector field along rays. 
The first one was not new, since it appeared in the very first paper about
Lienard equation \cite{L}, which gave to such an equation its name. The second one
was probably used for the first time just by Sansone, as well as the third one, even if this approach is
usually credited to Massera for his geometrical idea centered on the rotation of the vector field on rays. 

 In section 2 we give an account of the evolution of the theory starting from Levinson-Smith paper. Some of the more recent developments in this direction will be also presented. We know that not all the authors which deserve to be quoted will appear and we apologize for this. 
 
In section 3 we discuss some problems arising from the geometrical approach to uniqueness. We also present a new result in the light of the ideas of Sansone and Massera. Such a result can be applied to a class of planar differential systems, not equivalent to second order differential equations. This extends a recent result by Ciambellotti also based on Massera's geometrical approach \cite{Cia}. 

In section 4 we point out that the field still uncovered by existing uniqueness theorems contains 
strikingly simple, but non-trivial equations. Lins-De Melo-Pugh's conjecture, still open for $n \geq 3$, resists to researchers' efforts since 1976.

 \section{Uniqueness theorems}

Let us consider a Li\'enard equation in the phase plane,
\begin{equation} \label{phaseplane}
\dot x = y , \qquad \dot y = - g(x) - f(x) y  .
\end{equation} 
We shall also consider the usual system, known as Li\'enard plane,
\begin{equation} \label{lienardplane}
\dot x = y -F(x), \qquad \dot y = - g(x )  ,
\end{equation} 
where $F(x) = \int_0^x f(s)ds$. Throughout all of this paper we assume $f(x)$ to be continuous  and $g(x)$ to be locally lipschitzian on their common domain, so to have existence and uniqueness of solutions.  In the following we also write $G(x) = \int_0^x g(s)ds$, referring to 
$$
E(x,y) = G(x) + \frac {y^2}2
$$
as the energy function of both systems (\ref{phaseplane}) and (\ref{lienardplane}). Some of the basic properties and techniques concerned with limit cycles' existence and uniqueness appear already in Li\'enard's paper \cite{L}, which may be considered a milestone in this field, but after this paper, the first relevant result concerning the uniqueness of limit cycles is due to Levinson and Smith \cite{LS}. In particular, the authors presented the following results. We write $F(+\infty)$ for $ \int_0^{+ \infty} f(x) dx$. Similarly for $G$.

\begin{theorem} \label{LS1} Assume that $xg(x) >0$ for $x \neq 0$. If
\begin{itemize} 
\item there exist $\delta^- < 0 < \delta^+$ such that $G(\delta^-) = G(\delta^+)$;
\item $f(x) < 0$ in $ (\delta^-,\delta^+)$;
\item  $f(x) \geq 0$ in $ (-\infty,\delta^-] \cup [\delta^+,+\infty)$;
\item $G(\pm \infty) = F(+ \infty) = + \infty $;
\end{itemize}
then the system (\ref{phaseplane}) has exactly one limit cycle. 
\end{theorem}

Such a theorem is obtained as a special case of a more general one concerning the system
\begin{equation}  
\dot x = y , \qquad \dot y = - g(x) - f(x,y) y  .
\end{equation} 

Restricting to classical Li\'enard systems, the same paper contains also  a result strongly depending on symmetry properties of $f(x)$ and $g(x)$. Even if it is not explicitly stated as a corollary, we can summarize the containt of section 4 in \cite{LS} in the following statement.

\begin{theorem} \label{LS2} Assume $f(x)$ even, $g(x)$ odd, and $xg(x) >0$ for $x \neq 0$. 
If 
\begin{itemize} 
\item there exists $x_0 > 0$ such that $F(x) < 0$ in $ (0,x_0)$;
\item $F(x) > 0$ and increasing in $(x_0,+\infty)$;
\item $G(+ \infty) = F(+ \infty) = + \infty $;
\end{itemize}
then (\ref{lienardplane}) has exactly one limit cycle.
\end{theorem}

The proof did not introduce relevant novelties, being based on  Li\'enard's approach.   In fact, theorem (\ref{LS2}), under the condition $g(x) = x$, was proved by Li\'enard in his celebrated paper \cite{L}. Both Li\'enard's and Levinson-Smith's proofs are based on the fact that the energy's derivative along the solutions of (\ref{lienardplane}) is $\dot E(x,y) =  -F(x)g(x)$, so that the integral of $ -F(x)g(x)$ along a cycle is zero. Splitting the integral in several steps, where the integration is performed with respect to $x$ or $y$, and comparing analogous sub-integrals of distinct cycles,  one proves that  if $\gamma_1$ and $\gamma_2$ are two concentric cycles, with periods $T_1$ and $T_2$, 
$$
0 = \int_0^{T_1} E(\gamma_1(t)) dt \neq \int_0^{T_2} E(\gamma_2(t)) dt = 0,
$$ 
so obtaining a contradiction. 

A few years after Levinson and Smith, Sansone proved the following theorem, where the symmetry assumption on $F(x)$ was replaced by a weaker one \cite{San1}.

\begin{theorem} \label{San1} If $g(x) = x$, and
\begin{itemize} 
\item there exist $\delta^- < 0 < \delta^+$ such that $f(x) < 0$ for $x \in (\delta^- , \delta^+)$;
\item either $f(x) > 0$ in $(\delta^+,+\infty)$ or $f(x) > 0$ in $(-\infty,\delta^-)$;
\item there exists $\Delta > 0$ such that $F(\Delta) = F(- \Delta) = 0$;
\item $F(+\infty) = + \infty$, or $F(-\infty) = -  \infty$;
\end{itemize}
then the system (\ref{lienardplane})  has exactly one limit cycle, which is stable. 
\end{theorem}

In the same direction, for more general classes of systems, we can find the recent  results of Xiao  and Zhang Zhi-fen \cite{XZ1},   \cite{XZ2}, Carletti and Villari \cite{CV}, Sabatini and Villari \cite{SV}, Carletti \cite{C}. 

Next theorem was based on a different principle \cite{San1}.

\begin{theorem} \label{San2} If $g(x) = x$, and
\begin{itemize} 
\item there exists $\delta > 0$ such that $f(x) < 0$ for $x \in (-\delta , \delta)$,   $f(x) > 0$ in  $(-\infty,-\delta) \cup (\delta ,+\infty)$;
\item $F(+\infty) = + \infty$, or $F(-\infty) = -  \infty$;
\end{itemize}
then the system (\ref{phaseplane}) has exactly one limit cycle, which is stable. 
\end{theorem}

The above theorem, even if a special case of theorem \ref{LS1}, introduced a stability argument often re-used  by several authors in successive papers. Sansone proves that if a $T$-periodic  limit cycle $\gamma(t)$ exists, then it is attractive, since the divergence integral
$\int_0^T {\rm div}(\gamma(t)) dt$
is negative. Then, since two adjacent concentric limit cycles cannot be both attractive, the uniqueness follows. 
Such a theorem, even if not explicitly stated, proves also the cycle's hyperbolicity, that is an important feature in relation to perturbation problems.

An equally innovative result obtained by Sansone on limit cycle's uniqueness is virtually unknown to most researchers. It was  exposed in a talk, together some other results about Li\'enard equation, whose containt  appeared in \cite{San2}, section 2. We report here the uniqueness theorem for Li\'enard equation.

\begin{theorem}  \label{San3} If $g(x) = x$, and
\begin{itemize} 
\item there exist $\delta^- < 0 < \delta^+$ such that $f(x) < 0$ for $x \in (\delta^- , \delta^+)$;
\item either $f(x) > 0$ in $(\delta^+,+\infty)$ or $f(x) > 0$ in $(-\infty,\delta^-)$;
\item $f(\delta^-) = f(\delta^+) = 0$;
\item  $|f(x)| < 2$; 
\item  $f(x)$ is non-increasing in $(-\infty,\delta^-)$ and non-decreasing in $(\delta^+,+\infty)$ ;
\end{itemize}
then the system (\ref{lienardplane}) has exactly one limit cycle. 
\end{theorem}

The proof was based first on the trasformation  of the system (\ref{lienardplane}) into polar coordinates $(r,\theta)$, then on the study of the function
$$
\frac{d}{d \theta}  \ln r(\theta) =  \frac{ f(r \cos \theta) \sin^2 \theta}{1+ f(r \cos \theta) \sin \theta \cos \theta}.
$$
A generalization of Sansone's result is due to Conti, who used the map  
$$
(x,y) \mapsto (sign(x) \sqrt{G(x)}, y)
$$
known as Conti-Filippov transformation, to transform the system (\ref{phaseplane}) or (\ref{lienardplane}) into a new system with $g(x)$ replaced by $x$.  

The weakness of last Sansone's result is the assumption $|f(x)| < 2$ (which avoids the denominator's vanishing), clearly not satisfied by polynomials. In Sansone's approach it was not clear whether such an obstacle could be overcome.
Such an assumption was removed by Massera, when he was visiting Florence and working with Sansone. Actually, Massera's main contribution consisted in turning Sansone's analytical approach into a geometrical one. Massera observed that the monotonicity conditions on $f(x)$ imply that the vector field rotates clockwise,  as $r$ increases, along rays contained in  the half-plane $x > 0$, counter-clockwise,  as $r$ increases, along rays contained in the half-plane $x < 0$.  If $\gamma$ is star-shaped, then ever curve $\gamma^\kappa$ obtained from $\gamma$ by a $\kappa$-dilation (homothety) is as well star-shaped. The family $\gamma_k$ is a foliation of the punctured plane, and can be assimilated to the family of level curves of a Liapunov function for $\gamma$. The vector field rotation forces the orbits of (\ref{phaseplane}) to enter such curves, approaching $\gamma$. This in turn allows to show that every limit cycle is attractive, so preventing the co-existence of concentric limit cycles, by Sansone's attractivity argument.  On the other hand, if the limit cycle is not star-shaped, the curves $\gamma_k$ can intersect with each other, and for values of $\kappa$ close to $1$ they even intersect the limit cycle $\gamma$, so that proving $\gamma$'s attractivity is not immediate.  It is amazing to observe that in the original paper Massera deals with this point by just writing {\it Il faut remarquer que le raisonnement resterai valable m\`eme si $\gamma$ n'est pas \'etoil\'e par rapport \`a l'origine: dans ce cas $\gamma$ et $\gamma_k$ pourraient avoir des intersections.  }

It is worth noticing that Sansone's condition $|f(x)| < 2$ implies that for every non-trivial orbit the angular speed does not vanish. In particular, this occurs at cycles, which hence  have to be star-shaped.

It is also surprising to observe that the cycle's star-shaped-ness could have been easily proved by using the monotonicity conditions on $f(x)$, as shown in  \cite{ZZ}, p. 225. 
A more geometrical proof of the cycle's star-shaped-ness was developed by Villari, and appeared for the first time in \cite{CRV}. A second objection to Massera's geometrical approach is that, in general, the vector field could rotate up to reach a position opposite to the original one (tangent to $\gamma$) so that the successive position would take the orbits to leave the external $\gamma_k$'s. This cannot occur for second order equations, as will be shown in section 3.

We emphasize that the cycle's star-shaped-ness does not imply that the angular speed be of constant sign on every orbit. In fact, even for the simplest example of Li\'enard equation with a limit cycle, i.e. Van der Pol's equation, 
$$
\ddot x + \epsilon (x^2-1) \dot x + x = 0,
$$
infinitely many orbits passing through the second and fourth orthants change angular speed. We may consider the information on the cycle's  star-shaped-ness also as information about the cycle's location: it is contained  in the region $x^2-xyf(x)+y^2 > 0$. 
Other extensions were given by Carletti, Rosati and Villari \cite{CRV} and  Ciambellotti \cite{Cia}, but in both papers the geometrical idea of Massera is not really modified. 

At the end of this review, we cite a result that renews the divergence approach in a non-traditional way.  In some recent papers new classes of so-called {\it stability operators} were introduced, in order to study cycles' hyperbolicity. Such operators are just functions $\Psi(x,y)$ with the property that if $\gamma$ is a $T$-periodic cycle and
$$
\int_0^T \Psi(\gamma(t)) dt \neq 0,
$$
then $\gamma$ is hyperbolic, the stability character being given by the sign of the above integral. The divergence is clearly one of such operators. Other ones are the curvature of the orthogonal vector field, or the function
$$
\nu =  \frac{[V,W] \wedge V}{V\wedge W},  \qquad V\wedge W \neq 0,
$$
where $[V,W]$ is the Lie bracket of $V$ and $W$, and $V\wedge W $ is
the determinant of the matrix having $V$ and $W$ as rows   (\cite{FGG2}, \cite{GS}). The function $\nu$ has the following relevant property 
$$
\int_0^T (\mbox{div }V)(\gamma(t))\, dt = \int_0^T  \frac{[V,W] \wedge
V}{V\wedge W}\left( \gamma(t) \right)\, dt,
$$
In particular, if $\nu \leq 0$ in a domain (not identically vanishing on any cycle), then such a domain cannot contain two concentric limit cycles, since they both should be attractive. This is just the original argument by Sansone, applied replacing the divergence with $\nu$.  We emphasize that using $\nu$ allows also to prove the cycle's hyperbolicity, which is not a consequence of Sansone-Massera's approach. 

Finally, we observe that  Massera theorem is a special case of such a result, obtained by taking $W(x,y) = (x,y)$. In fact, as shown in Corollary 6 of \cite{GS}, given an arbitrary planar differential system
\begin{equation} \label{sistema}
\dot x = P(x,y), \qquad \dot y = Q(x,y)  ,
\end{equation} 
taking $W(x,y) = (x,y)$, gives, for $Q \neq 0$,
\begin{equation}\label{nuPQAB}
\nu =  \frac{P(xQ_x + yQ_y) - Q (xP_x + yP_y)}{y P - x Q} = \frac{r \left(  \dr \frac QP \right) }{y P - x Q} =\frac{r}{y \dot x - x \dot y} \left(  \dr \frac {\dot y}{\dot x} \right) . 
\end{equation}
In other words, if $\nu$ does not change sign, then the vector field rotates monotonically along rays as $r$ increases. 
Moreover, if  
$$\
\dot x  = P(x,y) = y, \qquad \dot y  = Q(x,y) = - x - y f(x) ,
$$
then
\begin{equation}  \label{massera}
 \nu =  - xf'(x)  \frac{y^2   }{x^2 + xyf(x) +y^2 } . 
\end{equation}
Hence, if the denominator $x^2 + xyf(x) +y^2$ does not vanish,  $\nu$'s sign is the same as that of $xf'(x)$, re-finding Sansone-Massera's condition.

 \section{A new result}

Let us call {\it radial angular monotonicity}, RAM, the property of the vector field to  rotate monotonically along rays as $r$ increases.
 An implicit assumption of Massera-like theorems is the fact that limit cycles rotate clockwise around the origin. This is a consequence of $\dot x =y$, which forces cycles cross the $y$-axis rotating clockwise.  If RAM holds, this gives the attractivity of every cycle, hence its uniqueness. In fact, the same assumption on a cycle rotating counter-clockwise generates repulsiveness, rather than attractiveness.  This allows as well to prove uniqueness, if RAM holds. On the other hand, the possible co-existence of  limit cycles rotating both clockwise and counter-clockwise does not allow to prove uniqueness under the only RAM hypothesis. Counterexamples with finitely many and infinitely many cycles have been given in \cite{S2}. We report here a polynomial system satisfying RAM, with two limit cycles rotating in different ways.
 \begin{equation}\label{sysduecicli} 
\left\{
\begin{array}{rl}
\dot x   = &\ \  y \Big(x^2+y^2-(x^2+y^2)^2 \Big)  + x \Big(1-3 (x^2+y^2)+(x^2+y^2)^2 \Big) \\ 
 \dot y  = & -x \Big(x^2+y^2-(x^2+y^2)^2 \Big) + y \Big(1-3 (x^2+y^2)+(x^2+y^2)^2 \Big)
\end{array}
\right.
\end{equation} 
Such a system has two star-shaped limit cycles coinciding with the circles $x^2+y^2 =\frac{3-\sqrt {5}}{2}$ and  $x^2+y^2 =\frac{3+\sqrt {5}}{2}$.  The internal one is an attractor, the external one is a repellor. The vector field rotates clockwise along every ray (see figure 1). 
\begin{figure}[h!]
  \caption{The system (\ref{sysduecicli}) has two limit cycles.}
  \centering
    \includegraphics[width=0.5\textwidth]{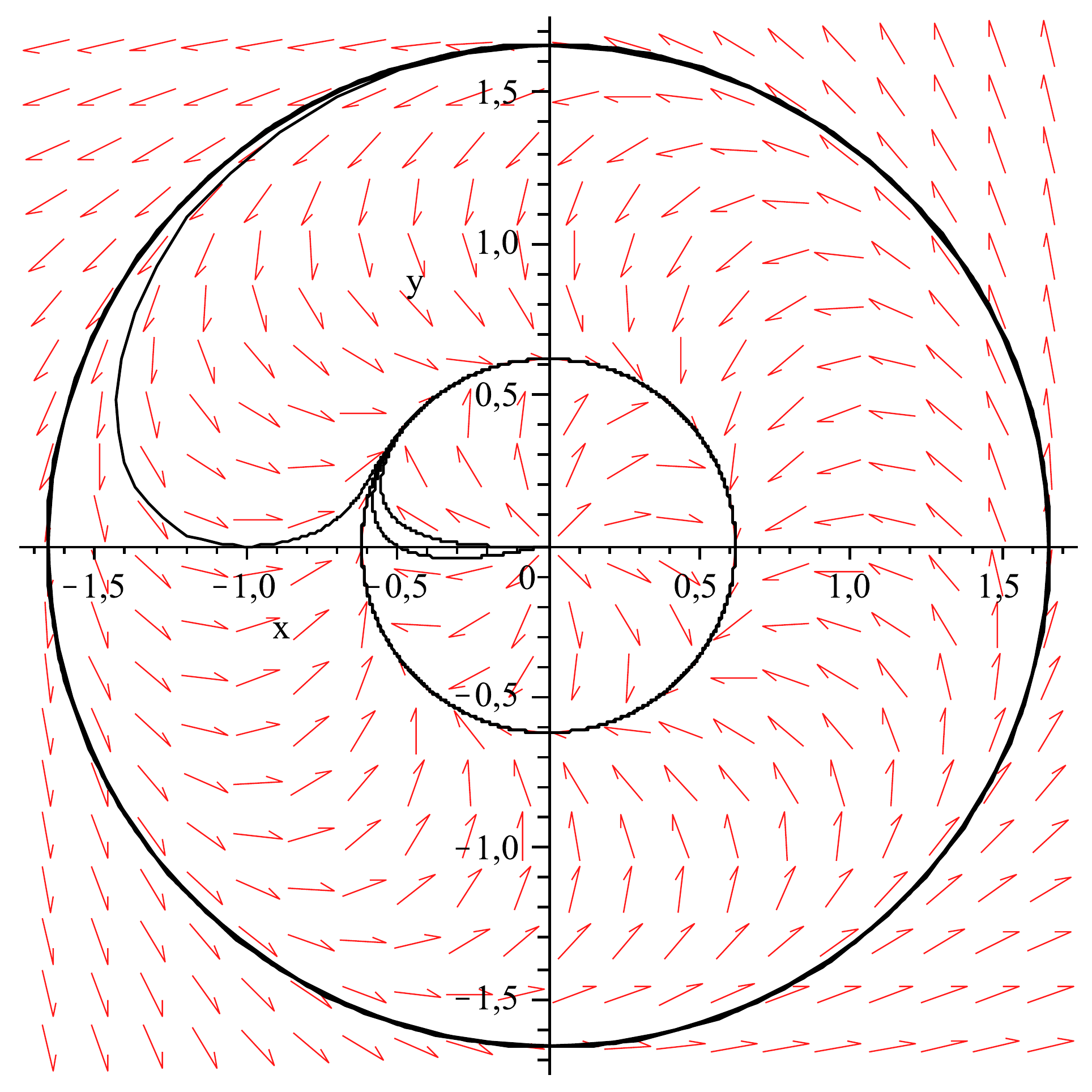}
\end{figure}

Hence, in order to prove a Massera-like theorem valid for systems, one has to assume some additional hypothesis that prevents the coexistence of opposite rotations. This does not necessarily have to hold in all of the plane. It is sufficient, to prove that it occurs in a region containing all cycles, or, as in  \cite{S1,S2}, in the region bounded by two adjacent limit cycles. In  \cite{S1,S2}, this was proved for systems of the type
$$
\dot x = y, \qquad \dot y = -x - y\phi(x,y)
$$
under the assumtpion that $\phi$'s level sets be star-shaped. This does not imply the orbits' star-shaped-ness.

We propose here some simple additional conditions that allow to prove uniqueness if RAM holds. We state next theorem for a system defined on the whole plane, but it can be easily adapted to an arbitrary star-shaped set containing the origin. We assume existence and uniqueness of solutions to (\ref{sistema}). Let us set 
\begin{equation}  
\alpha(x,y) = P( xQ_x  + yQ_y) - Q( xP_x  + yP_y).
\end{equation}
In next theorem, the interval $(a,b)$ may be a half-line or even all of $\R$.

\begin{theorem} \label{teorema}
Assume that  (\ref{sistema}) to have a unique equilibrium point at the origin. 
Assume $\alpha(x,y) \geq 0$ in $ \Omega $.  If one of the following conditions holds
 \begin{itemize}
 \item[(1)] for all $(x,y) \in \R^2 $, $(x,y) \neq (0,0)$, for all $\tau >0$, the vectors $V(\tau x,\tau y)$ and $V(x, y)$ are  linearly independent; 
 \item[(2)] there exists a $C^1$, open Jordan curve curve $\eta: (a,b) \rightarrow \R^2$ such that 
 $\lim_{s \rightarrow a^+}|\eta(s)| = 0$, $\lim_{s \rightarrow b^-}|\eta(s)| = +\infty$, and $\dot \theta(\eta(s)) > 0 \ \ \ (< 0)$;
 \end{itemize}
then  the system (\ref{sistema}) has at most one limit cycle in $\Omega$.
\end{theorem}
{\it Proof.}  According to (\ref{nuPQAB}), the sign of $\alpha$ is the sign of the radial derivative of $\displaystyle{\frac {\dot y}{\dot x} }$. Working as in \cite{CRV}, one proves that every cycle is star-shaped The existence of a single equilibrium point implies that all cycles are concentric. Then we work as follows.

Point (1). It is a modification of the original argument by Massera.  If a cycle $\gamma$ exists, as in \cite{M} we may consider the family of curves $\gamma_\kappa$ obtained from $\gamma$ by means of $\kappa$-dilations (homotheties). In particular, $\gamma = \gamma_1$. Without loss of generality, we may assum $\gamma_1$ to rotate clockwise. Then, the  vector field's rotation along rays implies that, for small positive values of $\kappa - 1$, $V$ points towards the interior of  $\gamma_\kappa$. Assume, by absurd, the existence of $\kappa^*$ and $(x^*,y^*)$ such that $V(x^*,y^*)$ points outwards $\gamma_{\kappa^*}$. Then, by continuity, there exists $\kappa^+$ and $(x^+,y^+)$ such that $V(x^+,y^+)$ is tangent to $\gamma_{\kappa^+}$. The rotation's monotonicity along rays   implies that this occurs first at a point $(x^+,y^+)$ where $V(x^+,y^+)$ is parallel and opposite to $V(x_1,y_1)$, with $(x_1,y_1) \in \gamma$. This contradicts the hypothesis (1). Internal attractivity can be proved in the same way. 

Point (2). Assume, by absurd, the existence of two distinct cycles  $\gamma_i$, $\gamma_e$, with $\gamma_i$ encircled by $\gamma_e$. Without loss of generality, we may assume them to be adjacent. Since both are star-shaped,  $\dot \theta$ does not change sign on any of them. Moreover, since they both cross the curve $\eta$, both have to contain a point where $\dot \theta > 0$. Hence $\dot \theta(\gamma_i) \geq 0$ and $\dot \theta(\gamma_e) \geq 0$. The condition RAM implies that both are attractive, which is a contradiction.
\hfill$\clubsuit$

Both hypotheses (1) and (2) are satisfied by Li\'enard systems. Point (1) since $\dot x = y$ on the curve $\eta(s) = (0,s)$, defined on $(0,+ \infty)$.

Next  corollary extend Ciambellotti's result \cite{Cia}. Consider the following family of systems
\begin{equation}   \label{homog}
\dot x = k(y), \qquad \dot y = -f(x) l(y) - \sum_{j=0}^n h_j(x) m_j(y).
\end{equation}

\begin{corollary} Assume  $k(y)$, $l(y)$,  to be $d$-homogeneous functions, $h_j(x)$  $j$-homogeneous, $m_j(y)$  $(d-j)$-homogeneous, for $j=1, \dots, n$, $f(x)$ differentiable. If $yk(y) >0$ and $yl(y) >0$ for $y \neq 0$, $xf'(x) \geq 0$ for $x \neq 0$, then the system   (\ref{homog}) has at most one limit cycle.
\end{corollary}
{\it Proof.}  The system  (\ref{homog}) satisfies the hypothesis (2) of  theorem  \ref{teorema}, taking $\eta(s) = (0,s)$, $s \in (0,+\infty)$.
Let us compute the function $\alpha$ for the system (\ref{homog}):
$$
\alpha=k\left[   -xf'l - x\sum h_j'm_j - y fl' - y\sum h_jm_j'\right] + \left[  f l + \sum_{j=0}^n h_j m_j  \right] yk' =
$$
$$
\ =  -xf' lk - k\sum xh_j'm_j -  fyl'k - k \sum h_jym_j'  +  f l yk' + yk' \sum_{j=0}^n h_j m_j  = 
$$
$$
\ =   -xf' lk - k\sum jh_j m_j -  dflk - k \sum h_j(d-j)m_j  + d f l k + dk \sum_{j=0}^n h_j m_j  
$$
$$
=   -xf' lk + k \left[  -\sum jh_j m_j - \sum h_j(d-j)m_j   + d\sum_{j=0}^n h_j m_j  \right] = -xf' lk .
$$
The sign of $-xf'(x) l(y)k(y)$ is the same as that of $-xf' (x)$, since both $k(y)$ and $l(y)$ have the sign of $y$. Hence the sign of $\alpha$ is that of $-xf'(x)$, as in formula  (\ref{massera}). 
\hfill$\clubsuit$

An example is given by the the system
\begin{equation}   \label{esempio}
\dot x = y^3, \qquad \dot y =  (5x^2-1) y^3 - x^3 -xy^2,
\end{equation}  
illustrated in figure 2.

\begin{figure}[h!]
  \caption{The system (\ref{esempio}) has a unique limit cycle.}
  \centering
    \includegraphics[width=0.5\textwidth]{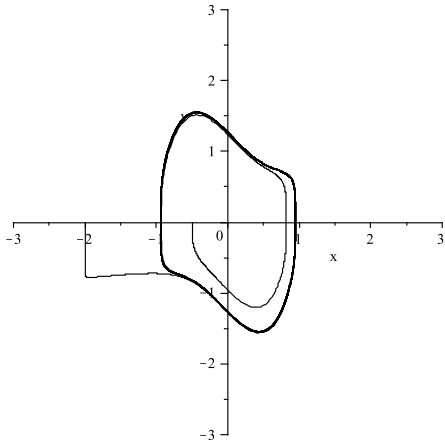}
\end{figure}

\newpage

 \section{An open problem}

At the end of this paper, after having listed several theorems of seemingly great generality, we would like to point out a very simple problem still waiting for a solution. Let us consider the Li\'enard equation with $f(x)$ cubic polynomial,
$$
\ddot x +(a x^3 + b x^2 + c x +d ) \dot x + x = 0, \qquad a,b,c,d \in \R.
$$
Such an equation is clearly out of the reach of any of the above results. Proving the limit cycle's uniqueness for such a class of equations would give an answer to Lins-De Melo-Pugh's Conjecture for $n = 3$ \cite{LDP}. As it is well known, such a conjecture is part of the famous XVI Hilbert's problem, re-proposed by Smale as one of the {\it Mathematical problems for the next century} \cite{Sm}.
The conjecture for $n=2$ was proved by Lins-De Melo-Pugh in their paper, but they were not aware of the fact that this problem was actually already solved by Zhang-Zhifen in 1958. This is a very deep result, but the for long time the proof was only in Chinese, until an English version appeared \cite{ZZ0}.

\end{document}